\documentclass{article}
\usepackage{xy}
\usepackage{epsfig}
\usepackage[T2A]{fontenc} 
\usepackage[cp1251]{inputenc} 
\usepackage[russian]{babel}
\usepackage{amsfonts,amssymb,amsmath}
\usepackage{graphics} 
\usepackage{multido, xcolor} 

\newtheorem{theorem}{\bf\large Теорема}[section]

\newtheorem{corollary}[theorem]{\large\bf Следствие}
\newtheorem{definition}[theorem]{\large\bf Определение}

\newtheorem{proposition}[theorem]{\large\bf Предложение}

\def\qed{ \ \hfil$\square$}

\def\Im{{\hskip0.3mm\rm  Im\hskip0.5mm}}

\makeindex

\def \F {{\mathbb F}}

\newcommand{\dd}{{\rm d}}
\newcommand{\cont}{{\rm cont}}

\newcommand{\GL}{{\rm GL}}

\newcommand{\Hom}{{\rm Hom}}

\newcommand{\SL}{{\rm SL}}

\newcommand{\tors}{{\rm tors}}

\newcommand{\N}{{\mathbb N}}
\newcommand{\Q}{{\mathbb Q}}
\newcommand{\R}{{\mathbb R}}
\newcommand{\Z}{{\mathbb Z}}

\newcommand{\Cr}{{\mathcal C}}

\newcommand{\Mc}{{\mathcal M}}

\newcommand{\Or}{{\mathcal O}}
\newcommand{\Pc}{{\mathcal P}}

\newcommand{\Sc}{{\mathcal S}}

\newcommand{\al}{\alpha}
\newcommand{\De}{\Delta}
\newcommand{\de}{\delta}
\newcommand{\e}{\varepsilon}

\newcommand{\Ga}{\Gamma}
\newcommand{\ga}{\gamma}

\newcommand{\La}{\Lambda}

\font\teneusm=eusm10 \font\seveneusm=eusm7 
\font\fiveeusm=eusm5 
\newfam\eusmfam 
\textfont\eusmfam=\teneusm 
\scriptfont\eusmfam=\seveneusm 
\scriptscriptfont\eusmfam=\fiveeusm

\def\mat #1,#2,#3,#4,{\left({#1\atop #3}{#2\atop #4}\right)}
\def\bra#1,{{\left\lbrace {#1}\right\rbrace}}

\def\si{\sigma}
\def\z{\zeta}

\def \Step {{\rm Step}}

\def\bs{\backslash}

\def\l1{\langle}

\newcommand{\B}{\left(\begin{array}{cc}}
\newcommand{\E}{\end{array}\right)}

\usepackage{euscript}
\let\scr=\EuScript

\let\mathcal=\scr           \def\Or{{\scr O}}
\def\ang#1,{{\left\langle {#1}\right\rangle}} 
\newcommand{\ds}{\displaystyle}

\def\qed{ \ \hfil$\square$}

\def\Im{{\hskip0.3mm\rm  Im\hskip0.5mm}}

\def \F {{\mathbb F}}

\def\fraco #1,#2,{{{#1}\over {#2}}}

\font\teneusm=eusm10 \font\seveneusm=eusm7 
\font\fiveeusm=eusm5 
\newfam\eusmfam 
\textfont\eusmfam=\teneusm 
\scriptfont\eusmfam=\seveneusm 
\scriptscriptfont\eusmfam=\fiveeusm

\font\tengothic=eufm10
\font\sevengothic=eufm7
\font\fivegothic=eufm5
\newfam\Gothic
\textfont\Gothic=\tengothic
\scriptfont\Gothic=\sevengothic
\scriptscriptfont\Gothic=\fivegothic
\def\go{\fam\Gothic\tengothic}

\def\mat #1,#2,#3,#4,{\left({#1\atop #3}{#2\atop #4}\right)}
\def\bra#1,{{\left\lbrace {#1}\right\rbrace}}

\def\si{\sigma}
\def\z{\zeta}
\def\w{\omega}

\def \Re {{\rm Re}}

\def \Card {{\rm Card}\,}
\def \Step {{\rm Step}}

\def\bs{\backslash}

\def\l1{\langle}
\def\Eq{\Longleftrightarrow}

\let\scr=\EuScript

\let\mathcal=\scr           \def\Or{{\scr O}}
\let\db=\mathbb

\def\dbQ{{\db Q}}

   \def\dbZ{{\db Z}}

\def\vin{{ {\tiny \mid }  
\kern-7.29pt 
\bigcup }}

\def \bs{{\backslash}}

\def\ang#1,{{\left\langle {#1}\right\rangle}}

\normalsize

\newcommand{\CC}{\mathbb C}

\newcommand{\GG}{\mathbb G}
\newcommand{\HH}{\mathbb H}

\newcommand{\QQ}{\mathbb Q}

\newcommand{\ZZ}{\mathbb Z}

\def\oc{\Or}

\newcounter{ncours}{\setcounter{ncours} {1}}

\def\oc{{\cal O}}

\def\SL {\mathop{\rm SL}\nolimits}

\font\tenrus=wncyr10
\font\sevenrus=wncyr7
\newfam\Rus
\textfont\Rus=\tenrus
\scriptfont\Rus=\sevenrus

\def\scr{\fam\Rus\sevenrus} 

\font\tenbrus=wncyb10
\font\eightbrus=wncyb8
\newfam\Brus
\textfont\Brus=\tenbrus
\scriptfont\Brus=\eightbrus

\font\tenirus=wncyi10
\font\eightirus=wncyi8
\newfam\Irus
\textfont\Irus=\tenirus
\scriptfont\Irus=\eightirus

\font\tenrus=wncyr10
\font\sevenrus=wncyr7
\newfam\Rus
\textfont\Rus=\tenrus
\scriptfont\Rus=\sevenrus

\def\scr{\fam\Rus\sevenrus} 
\input cyracc.def
\def\i1{\accent'044i}
\def\I1{\accent'044I}
\def\e1{\accent'040e}
\def\l1{l{}p1}

\title{Модулярные формы и $p$-адические числа
}

\author 
{А.А.Панчишкин
}
\date{}

\def\ZZ{\mathbb Z}

\def\qed{\quad\hbox{\hskip 1pt\vrule width 4pt height 6pt
          depth 1.5pt\hskip 1pt}}

\begin{document}

\maketitle
%
%
\begin{abstract}
Пусть $p$ -- простое число. 
Обсуждаются  $p$-адические свойства различных арифметических функций, связанных с коэффициентами модулярных форм и производяшими  функциями.
Модулярные формы рассматриваются как средство решения задач арифметики.
Приведены примеры сравнений между модулярными формами, а также  
примеры компьютерных вы\-чис\-ле\-ний с модулярными формами и $p$-адическими числами.
\end{abstract}
  \tableofcontents
\section{Введение}

Статья основана на материалах спецкурсов автора в Университете Жозеф Фурье (Гренобль, Франция), 
 лекций автора  для французских педагогов в Институте Исследований Математичес\-кого Просвещения (IREM, Гренобль, Франция) в 1998, в 
Эколь Нормаль (Лион, Франция),
а также на материалах спецкурсов на мех-мате МГУ в 1979-1991 и в 2001.

В статье обсуждаются следующие темы:

\begin{enumerate}
\item[1)] Примеры производящих функций, 
модулярные формы и сравнения.
Представление целых чисел квадратичными формами.
 \item[2)] 
  Ряды Эйзенштейна и сравнения для функции Рамануджана.

\item[3)] 
{Числа и многочлены Бернулли, сравнения Куммера} 

\item[4)]
{Мера Мазура  и $p$-адическое интегрирование.}
\end{enumerate}
\section{Производящие функции, модулярные формы и сравнения.}
\subsection{Производящие функции}
Традиционной областью применения производящих функций является комбинаторика и теория разбиений. 
Пусть $p(n)$ — число разбиений
натурального числа $n$ в сумму натуральных неубывающих слагаемых:
\begin{align*}
&p(1)=1\quad :\quad  1=1;\cr
&p(2)=2\quad : \quad 2=2, \quad 1+1;\cr
&p(3)=3\quad : \quad 3=3,\quad 2+1,\quad 1+1+1;\cr
&p(4)=5;\quad p(5)=7.
\end{align*}
Тогда для производящей функции для $p(n)$ справедливо тождество Эйлера:
\begin{align}
1+\sum_{n=1}^\infty p(n)q^n=\prod_{m=1}^\infty (1-q^m)^{-1}. \label{i12.27}
\end{align}
Действительно, непосредственное перемножение показывает, что
\begin{align*}&
\prod_{m=1}^\infty (1-q^m)^{-1}= \prod_{m=1}^\infty (1+q^m+q^{2m}+q^{3m}+\cdots )= \\ &
(1+q+q^{2}+q^{3}+\cdots )\times (1+q^2+q^{4}+q^{6}+\cdots )\times \cdots \\&
\cdots\times (1+q^k+q^{2k}+q^{3k}+\cdots )\times \cdots 
=\sum_{a_1\ge 0, a_2\ge 0, a_3\ge 0, \cdots} q^{a_1+2a_2+3a_3+\cdots},
\end{align*}
а $p(n)$ как раз и есть число решений целых числах $a_1, a_2, a_3, ... , > 0$
«уравнения с бесконечным числом переменных» 
$$
a_1+2a_2+3a_3+\dots = n.
$$
Оказывается,  что бесконечные произведения типа (\ref{i12.27}) тесно связаны с тэта-рядами. 
 Например, при
$|q|<1,\ {u}\neq 0$ имеем(см. \cite{Andrews G.E. (1976)})
$$
\sum_{n=-\infty}^\infty {u}^nq^{n^{2}}=\prod_{m=0}^\infty (1-q^{2m+2})
(1+{u}q^{2m+1})(1+{u}^{-1}q^{2m+1})\qquad (\mbox{Якоби}), 
$$
$$
\sum_{n=0}^\infty q^{n(n+1)/2}={\ds\prod_{m=1}^\infty 
(1-q^{2m})\over \ds\prod_{m=1}^\infty (1-q^{2m-1})} 
\qquad (\mbox{Гаусс}),
$$
которые  выводятся из более общего тождества Коши:
при $|q|<1,\ |t|<1,\ a\in {\CC}$:
\begin{align}
1+\sum_{n=1}^\infty \frac{(1-a)(1-qa)\dots (1-aq^{n-1})t^n} 
{(1-q)(1-q^2)\dots (1-q^n)}  = \frac{\ds\prod_{m=0}^\infty(1-atq^m)} 
{\ds\prod_{m=0}^\infty(1-tq^m)} . \label{i12.28}
\end{align}
Вот иллюстрация вычисления с PARI-GP
(см. \cite{BBBCO}):
\begin{verbatim}
gp > {n=100;\\
prod(i=1,n,(1-x^(2*i)))*prod(i=1,n,((1-x^(2*i-1))^(-1))+O(x^(n+1)))
}
%5 = 1 + x + x^3 + x^6 + x^10 + x^15 + x^21 + x^28 + x^36 + x^45 + x^55 + x^66 +
 x^78 + x^91 + O(x^101)
gp > ##
  ***   last result computed in 451 ms.
  \end{verbatim}
\subsection{Представление целых чисел 
квадратичными формами.
}\label{i12.6.}
Пусть
\begin{align*}
&f(x)=f(x_1,\dots ,x_n)=\sum_{i,j=1}^{n}a_{ij}x_ix_j=A[x]=x^tAx,\cr
&g(y)=g(y_1,\dots ,y_m)=\sum_{i,j=1}^{m}b_{ij}y_iy_j=B[y]=y^tBy,
\end{align*}
—целочисленные
квадратичные формы с матрицами 
$A$ и $B$. 
Будем говорить, что квадратичная форма $f$ {\it 
представляет
$g$ над} ${\Bbb Z}$ если для некоторой целочисленной матрицы $C\in M_{n,m}(\Z )$ выполнено тождество
\begin{align}
f(Cy)=g(y), \quad {что эквивалентно}\quad A[C]=B. \label{i12.20}
\end{align}
В частности, при $m=1$ и $g(y)=by^2$, $f$ представляет форму $g$ если
$f(c_1,\dots ,c_n)=b$ для некоторых целых
чисел $c_1,\dots ,c_n$.

Лагранж доказал, что всякое целое число представимо суммой четырёх квадратов.
Этот факт выводится также из (более трудной) теоремы Гаусса о том, что целое положительное число $b>0$ 
тогда и только
тогда является суммой трех квадратов, когда оно не является числом
вида $4^k(8l-1),\quad
k,l\in \Z $
(см. \cite{Serre J.--P. (1970)}, \cite{Ma-Pa05}).

Пусть
\begin{align}
r_k(n)=\Card \{(n_1,\dots ,n_k)\in \Z ^k\ |\ n_1^2+\dots +n_k^2=n\}.
\label{i12.21}
\end{align}
число представлений $n$ в виде суммы $k$ квадратов. Так, например, $r_2(5)=8$, поскольку
\begin{align*}&
5 = 2^2 + 1^2 = 2^2 + (-1)^2 = (-2)^2 + 1^2 = (-2)^2 + (-1)^2 =\\ &
= 1^2 + 2^2 = (-1)^2 + 2^2 = 1^2 + (-2)^2 = (-1)^2 + (-2)^2. 
\end{align*}
В большом числе случаев найдены формулы для чисел представлений.
Приведем
лишь классический результат Якоби,
(см. \cite{Andrews G.E. (1976)}, \cite{Ma-Pa05}):
\begin{align}
r_4(n) = 
\begin{cases} \ds 8  \sum_{d|n}d,\quad &\hbox{если} \quad n 
\quad \hbox{нечётно},\cr \ds
24 \sum_{\scriptstyle d|n\atop \scriptstyle d\equiv 1(2)}d,\quad &\hbox{если} \quad n \quad \hbox{чётно}.
\end{cases} \label{i12.22}
\end{align}
 из которого также следует теорема Лагранжа.
Метод доказательства этой теоремы основан на введении производящей
функции для чисел $r_k(n)$:
$$
\sum_{n=0}^{\infty}r_k(n)q^n=\sum_{(n_1,\dots ,n_k)\in \Z^k }
q^{n_1^2+n_2^2+\dots + n_k^2}=\theta ({z})^k,
$$
где
\begin{align}
\theta ({z})= \sum_{n\in \Z }q^{n^2},\qquad q=e^{2\pi i{z}}.\label{i12.23}
\end{align}
— тэта-функция, которая рассматривается как голоморфная функция
на верхней комплексной полуплоскости
 $\HH=\{{z} \in {\CC}\ | \ \Im({z} )>0\}$
и обладает рядом замечательных аналитических свойств.
Эти свойства позволяют однозначно охарактеризовать $\theta^4({z} )$ как  {\it 
модулярную форму веса} 2 относительно группы $\Gamma_0(4)$, где используется обозначения  
\begin{align}
 \Gamma_0(N)
 =
 \left\{ \left( \begin{matrix}a & b\cr c & d \end{matrix}\right) \
 \in \ SL(2, \Z )\ \Bigg| \ N|c \right\}\subset \SL(2,{\Z}).\label{i12.24}
\end{align}
Другими словами, голоморфный дифференциал $\theta^4({z} )d{z}$ 
не меняется при дробно-линейных преобразованиях ${z} \mapsto (a{z} +b)
(c{z} +d)^{-1}$ с матрицей $\left(\begin{matrix}a & b\cr c & d 
\end{matrix}\right)$ in $\Gamma_0(4)$
(и  удовлетворяет оценкам регулярности роста при $\Im({z}) \to \infty$ в вершинах; заметим, что $2\pi i d{z}= {dq\over q}$, позтому дифференциал мероморфен с простым полюсом в точке $q=0 \Eq {z}=i\infty$).  


\subsection{Мотивировка: функция Рамануджана  $\tau$ и её контекст.} 
 
\medskip\medskip 
В качестве иллюстрации к общей теории приведём 
несколько удивительных свойств функции Рамануджана  $\tau$. 
 
\medskip 
Этот знаменитый пример происходит из происходит из следуюшей производяшей функции, 
определённой разложением в ряд следующего бесконечного произведения:
$$
q\prod_{m\ge 1}{(1-q^m)}^{24}=\sum_{n\ge 1}\tau
(n)q^n=q-24q^2+252q^3-1472q^4+\cdots
$$
 
\smallskip\noindent 
Положим $q={\rm exp}(2i\pi z)$ для $z$ 
из верхней комплексной полуплоскости
${{\HH}}=\{z\in
{{\CC}}\mid {\rm Im}(z)>0\}$, это
голоморфное отображение ${{\HH}}$
 на единичный круг с выколотым центром 
$q:{{\HH}}\rightarrow D(0,1)\backslash
\{0\}$. 
\par 

\medskip
Определяется функция  $\Delta :{{\HH}}\rightarrow {{\CC}}$, 
голоморфная на ${{\HH}}$, по формуле:   
\medskip
$$
\Delta (z)= \Delta _\infty (q)=q\prod_{m\ge 1}(1-q^m)^{24}
$$
 
\smallskip\noindent 
Эта функция даёт пример модулярной формы.  
Она обладает рядом замечательных свойств:

\subsubsection*{Автоморфность} \label{ens1.2.1} 
Группа $\SL(2,{\Z})$ целочисленных  квадратных $2\times 2$-матриц с определителем 1
 действует на ${{\HH}}=\{z\in
{{\CC}}\mid {\rm Im}(z)>0\}$
 по формуле
$$
\forall\gamma=\mat a,b,c,d,\in {\rm SL}(2,{\Z}),\ \ \biggr(\mat a,b,c,d,,z\biggr)\mapsto
\displaystyle \gamma \cdot z ={{az+b}\over{cz+d}}.
$$ 
\medskip 
Свойство автоморфности имеет вид:
\medskip\noindent
\begin{align}\label{Aut}
\forall\gamma=\mat a,b,c,d,\in {\rm SL}(2,{\Z}),\
\forall z\in {{\HH}} \Rightarrow \Delta (\gamma\cdot z)=(cz+d)^{12}\Delta (z).
\end{align}
\bigskip\noindent
Заметим, что свойство автоморфности (\ref{Aut})равносильно тому, что, голоморфный дифференциал $\Delta ({z})(d{z})^6$ 
не меняется при дробно-линейных преобразованиях ${z} \mapsto (a{z} +b)
(c{z} +d)^{-1}$ с матрицей $\left(\begin{matrix}a & b\cr c & d 
\end{matrix}\right)$ in ${\rm SL(2,{\Z})}$, поскольку для всех 
$\gamma \in {\rm SL(2,{\Z})}$, и для всех  ${z}\in {{\HH}}$, имеем   \par 
\medskip\noindent
$$
\gamma=\mat a,b,c,d,\Rightarrow d(\gamma\cdot
{z})=(c{z}+d)^{-2}d{z}.
$$
Отсюда непосредственно вытекает, что для любого натурального $m$ группа ${\rm SL(2,{\Z})}$ действует на множестве голоморфных функций $f(z)$на ${z}\in {{\HH}}$ 
по формуле: 
для  $\gamma \in {\rm SL(2,{\Z})}$, и для ${z}\in {{\HH}}$, имеем   \par 
\medskip\noindent
$$
\gamma=\mat a,b,c,d,\Rightarrow ( f|_{2m}\gamma)({z})=(c{z}+d)^{-2m}f(\gamma\cdot {z}),
$$
(действие веса $2m$),
а свойство автоморфности (\ref{Aut}) означает, что $\forall \gamma=\mat a,b,c,d,\Rightarrow \De|_{12}\gamma=\De$.
Поэтому (\ref{Aut}) достаточно проверить на образующих группы ${\rm SL(2,{\Z})}$.
%
Используем тот факт, что группа
$SL(2,{\Z})$ порождена матрицами 
$T= \mat 1,1,0,1,$ и  $S=\mat 0,-1,1,0,$. 
Чтобы в этом убедиться, используется 
алгоритм Евклида применительно к паре
 $(a,b)$, а также степени элемента
$S$, имеющего порядок  $4$.

Отсюда выводится, что свойство автоморфности (\ref{Aut}) достаточно проверять для элементов $S$ и  $T$, т.е. 
$$
\Delta (z+1)=\Delta (z),\ \ \Delta (-1/ z)=z^{12}\Delta (z),
$$
см. ниже.

\subsubsection*{Мультипликативность.} \label{ens1.2.2} 
Функция Рамануджана $\tau$ мультипликативна в следуюшем смысле
[обозначим через {\bf P} множество всех простых чисел]:  
\medskip\noindent
$$
\begin{cases}\forall m\in {\N}^*,\ \forall n\in {\N}^*,\ (m,n)=1
\Rightarrow \tau
(mn)=\tau(m)\cdot  \tau (n); \cr \forall p\in {\bf P},\ \forall r\in {\N}^*,\
\tau (p^{r+1})=
\tau (p^r)\tau (p)- p^{11}\tau (p^{r-1}); \cr \forall m\in {\bf
N}^*,\ 
\forall n\in {\N}^*,\ \tau (m) \tau (n)= \displaystyle\sum_{d\mid (m,n)}
d^{11}\tau(mn/d^2).
\end{cases}
$$
\smallskip\noindent
Эти свойства были предположены Рамануджаном и доказаны Морделлом и Гекке. 
Возможно, однако, что не существует ``элементарного''  доказательства этих свойств,  
в духе теоремы Гёделя о недоказуемости средствами элементарной арифметики, см.\cite{Ma-Pa05}. 
Может оказаться, что же замечание относится и к теореме Ферма, доказаной Уайлсом в 1994 
в высшей степени ``неэлементарными'' методами
[включающими теорию модулярных форм,
 $p$-адический анализ,
 теорию деформаций представлений Галуа,
 алгебраическую геометрию, 
   ...].

\medskip 
Естественная формулировка свойств мультипликативности функции Рамануджана
использует ряд Дирихле, связанный с функцией  $\tau$:
\medskip\noindent
$$
L(\Delta ,s)=\sum_{n\ge 1} \tau (n)n^{-s}=\prod_{p\in {\bf P}} \bigr(1-\tau
(p)p^{-s}+p^{11-2s}\bigr)^{-1}.
$$
\smallskip\noindent 
Этот ряд аналогичен ряду Дирихле задающему дзета-функцию Римана,
\medskip\noindent 
$$
\zeta (s)=\sum_{n\ge 1} n^{-s}=\prod_{p\in {\bf P}} \bigr(1-p^{-s}\bigr)^{-1},
$$
\smallskip\noindent
где равенство выражает свойство существования и единственности разложения 
натурального числа в произведение простых чисел.

Точно так же и в случае функции  Рамануджана $\tau$, справедливо тождество 
$$
\sum_{n\ge 1} \tau (n)n^{-s}=\prod_{p\in {\bf P}} \Biggr(\sum_{r\ge 0}
\tau\bigr(p^r\bigr)p^{-rs} \Biggr),
$$ 
а доказательство рекуррентных формул сводится к равенству: 
$$
\bigr(1-\tau (p)p^{-s}+p^{11-2s}\bigr)\cdot\Biggr(\sum_{r\ge 0} \tau
(p^r)p^{-rs}\Biggr)=1
$$ 

\subsubsection*{Оценки.}\label{ens1.2.3} 
 
Следующее свойство, первоначально предположенное Рамануджаном, 
 было доказано Делинем:
$$
\forall p\in {\bf P},\ \vert\tau (p)\vert <2p^{11/2}.
$$ 
Это свойство эквивалентно отрицательности дискриминанта многочлена
второй степени
 $X^2-\tau (p)X+p^{11}$, 
для всех простых чисел $p$. 
Для фиксированного $p$, пусть
$\alpha _p$ и $\beta _p$ -- комплексно-сопряжённые корни этого многочлена. 
Из формулы мультипликативности следует, что:  
\medskip\noindent
$$
{1\over{(1-\tau (p)X+p^{11}X^2)}}={1\over{(1-\alpha _pX)(1-\beta _pX)}}=
\Biggr(\sum_{r\ge 0} \tau (p^r)X^r\Biggr).
$$
Для всех  $r\ge 1$ выводится соотношение
$\tau (p^r)=\displaystyle\sum_{j=0}^r \alpha _p^j 
\beta _p^{r-j}=\sum_{j=0}^r \alpha _p^{2j-r} p^{11(r-j)}$. 
Абсолютная величина  $\alpha _p$ равна $p^{11/2}$, откуда следует оценка
$$
\vert\tau (p^r)\vert <(r+1)p^{11r/2}.
$$
Применение формального тождества даёт следующую оценку 
$$
\forall n\in {\N}^*,\ \vert\tau (n)\vert <\sigma _0(n)n^{11/2}=O(n^{{11
\over
2}+\varepsilon})
$$
\smallskip\noindent 
где $\sigma _0(n)$ -- число делителей числа $n$, где
$O({\rm ln}(n))=O(n^{\varepsilon})$ для любого $\varepsilon >0$.

Отсюда, в частности, выводится, что ряд $L(\Delta ,s)$ абсолютно сходится и определяет голоморфную функцию в правой полуплоскости Re$(s)>13/2$. 

\subsubsection*{Функциональное уравнение
для  $L(\Delta ,s)$.}\label{ens1.2.4} 
Определим функцию $L^*(\Delta ,s)$ по формуле
 $L^*(\Delta ,s)=(2\pi )^{-s}\Gamma (s)L(\Delta ,s)$. 
Эта функция, с одной стороны, продолжается до голоморфной функции 
на всей комплексной плоскости ${{\CC}}$, и эта функция удовлетворяет 
функциональному уравнению $L^*(\Delta ,12-s)=L^*(\Delta ,s)$. 
Можно сравнить это функциональное уравнение с функциональнум уравнением для дзета-функции Римана 
$\zeta (s)$
$$
\zeta ^*(s)=\pi ^{-s/2}\Gamma (s/2)\zeta (s)=\zeta^*(1-s).
$$
 
\subsubsection*{Связь с числами разбиений.}\label{ens1.2.5}  
Напомним,что разбиением натурального числа  $n$ называется неубывающая последовательность натуральных чисел с суммой, равной $n$. 
Функция числа разбиений обозначается через  $p:{\N}\rightarrow 
{\N}$ причём полагают $p(0)=1$. 

Как мы видели, производящий ряд функции
$p:{\N}\rightarrow {\N}$ даётся бесконечным произведением 
$\displaystyle \sum_{n\ge
0} p(n)q^n=\prod_{m\ge 1} (1-q^m)^{-1}$. 
Соответствующая голоморфная функция переменной $q$ сходится в открытом круге. 
Поэтому получается  голоморфная функция на ${{\HH}}$, $f:{{\HH}}
\rightarrow {{\CC}}$ где
$$
f(q)=\sum_{n\ge 0} p(n)q^n=\prod_{m\ge 1} (1-q^m)^{-1}.
$$
Имеет место равенство $\tilde\Delta (q)=q(f(q))^{-24}$ связываюшее функцию числа разбиений и функцию Рамануджана $\tau$. 

Используя свойство автоморфности, Харди и Рамануджан доказали следующую оценку для  $p(n)$:
$$
p(n)=\biggr({1\over{4\sqrt{3}}}+O\biggr({1\over{\lambda (n)}}\biggr)\biggr)\cdot
{{{\rm exp} (K\cdot\lambda (n))}\over{\lambda (n)^2}}
$$
где $\lambda (n)=\sqrt{n-{1\over 2}}$ и $K=\pi\sqrt{2/3}$ (см. \cite{Chand70}).
  
\subsubsection*{Сравнение Рамануджана и представления групп Галуа.} \label{ens1.2.6} 
Сравнение Рамануджана утверждает, что 
\begin{align}\label{Ram691}
\forall n\in {\Z}^+,\ \tau (n)\equiv\sum_{d\mid n} d^{11} {\bmod}\ 691.
\end{align}
В частности, $\tau(691) = -2747313442193908 \equiv 1 {\bmod}\ 691$.

Достаточно проверить  справедливость  сравнения $\tau (p)\equiv 1+p^{11}{\rm mod}\ 691$
для любого простого числа  $p$ отличного от  $691$. 
Действительно, тогда в силу мультипликативности и в силу рекуррентных соотношений по $r$
будем иметь
$$
\tau (p^{r+1})\equiv \tau (p^r)\tau (p) - p^{11}\tau (p^{r-1})
\equiv\sum_{j=0}^{r+1}
p^{11j} {\bmod} \ 691,
$$
откуда будет следовать и общее сравнение (\ref{Ram691}).
Серр нашёл объяснение этого курьёзного сравнения в рамках теории представлений групп Галуа.

\medskip
Пусть $\overline {\Q}$ -- алгебраическое замыкание поля   ${\Q}$ рациональных чисел. 
Пусть  $p$ -- простое число, отличное от  $691$ и ${\go p}$  --
произвольный простой идеал над $(p)$ 
в кольце ${\cal O}$ целых элементов $\bar {\Q}$. 
Обозначим через $G_{\go p}$ и $I_{\go p}$ 
 подгруппы группы Галуа
$G={\rm Gal}(\bar {\Q}/{\Q})$ определенные равенствами:
\begin{align*}&
G_{\go p}=\{\sigma\in G\mid\sigma {\go p}=
{\go p}\}\\ &
I_{\go
p}=\{\sigma\in G_{\go p} \mid\forall  x\in {\cal O},\ \sigma x\equiv x\ {\bmod}\ {\go
p}\}.
\end{align*}

\medskip
Группа $G_{\go p}$ называемая группой разложения,
 отождествляется с группой Галуа алгебраического замыкыния
 $\overline {\Q}_p$ поля ${\Q}_p$  $p$-адических чисел, 
 её нормальная подгруппа  $I_{\go p}$ называется группой инерции,
 а фактор-группа  
 $G_{\go p}/I_{\go p}$  отождествляется с группой Галуа
 алгебраического замыкыния конечного поля ${\F}_p$.
Эта группа порождается элементом Фробениуса  $Fr_p$.

\medskip 
Серр предположил, а Делинь доказал, что для любого простого числа $l$, существует такое представление Галуа
 $\rho _l:G\rightarrow {\rm GL}(2,{\Z}_l)$, что
для любого простого числа  $p$ отличного от  $l$, группа инерции  $I_{\go p}$ тривиально действует 
 (т.е. $\rho _l$ неразветвлено в  $p$), и  
det(Id$-\rho_l(Fr_p)\cdot X)=1-\tau (p)X+p^{11}X^2$. 
В случае $l=691$ справедливо сравнение
  $\rho _l(Fr_p)\equiv \mat p^{11} , \star ,  0 , 1, {\bmod}\ 691$, 
откуда $\tau (p)\equiv
1+p^{11} {\rm mod}\ 691$.

Впоследствии,  такие представления Галуа послужили основой для доказательства
теоремы Уайлса о модулярности эллиптических кривых (1994),
а также гипотез Серра о модулярности всех нечётных двумерных представлений Галуа над конечными полями, доказанных в 2007 Ч.Кхаре и Ж.-П. Вин\-тен\-бер\-же с использованием методов М.Кисина, Ж.-М. Фонтэна и Р.Тэй\-лора (Летняя школа в Марселе-Люмини, июль 2007).   
\begin{figure}[h]
\begin{center}
\includegraphics[width=14cm,height=10cm]{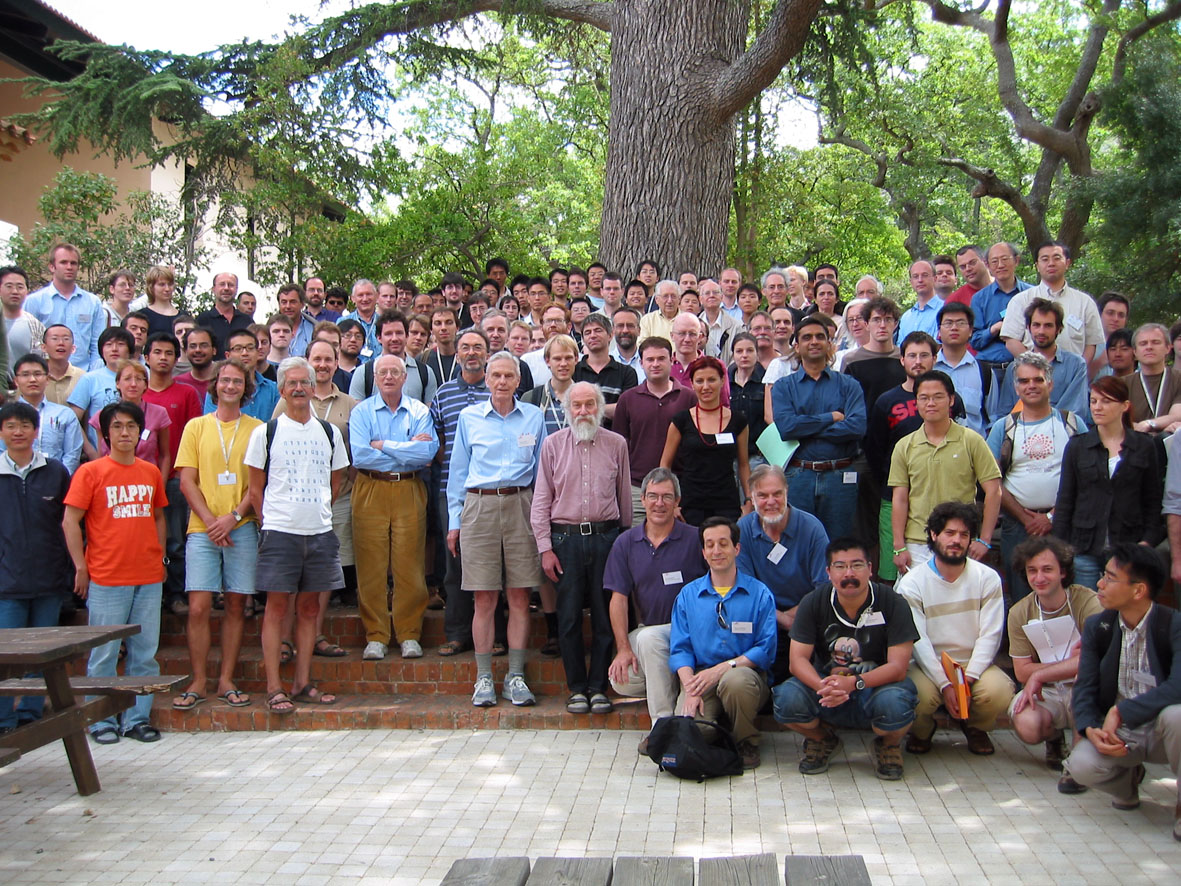}
\caption{\label{23w1_2}Летняя школа ``Гипотезы Серра о модулярности'' в Марселе-Люмини, июль 2007} 
\end{center}
\end{figure}

\subsubsection*{Формулы Ю.И.Манина}\label{ens1.2.7} 
Используя цепные дроби и модулярные символы, Ю.И.Манин
нашёл  формулы для функции Рамануджана 
 $\tau (n)$, даюшие гораздо более быстрый метод вычисления этой функции,
 чем метод разложения в ряд бесконечного произведения, или же метод основанный на рядах Эйзенштейна
 ($\De=(E_4^3-E_6^2)/1728$).     
 Эти формулы таковы:
\medskip\noindent
$$
\tau (n)=\sigma _{11}(n)-
 \sum {}^{{}^*(n)} 
\biggr({{691}\over{18}} 
(\Delta ^8\delta ^2-\Delta ^2\delta
^8)-{{691}\over{6}} (\Delta ^6\delta ^4-\Delta ^4\delta ^6)\biggr);
$$
\medskip\noindent
$$
\tau (n)=\sigma _{11}(n)-{{691}\over{18}}
 \sum {}^{{}^*(n)} 
\Delta ^2\delta ^2(\Delta^2-\delta ^2)^3
$$  
где $\sigma _{11}(n)=\displaystyle\sum_{d\mid n} d^{11}$ 
а во внешней сумме $ \sum{}^{{}^*(n)} $ справа суммирование производится по всем целым решениям уравнения,
$n=\Delta\Delta '+\delta\delta '$, 
котрые ``допустимы'', т.е. удовлетворяют условиям 
\begin{align*}
\big \{ (\De, \de)\big| &n=\Delta\Delta '+\delta\delta ',\ 
\De>\de>0,\ \De^\prime>\de^\prime>0, \ 
\hbox{где} \cr &
 \De | n, \ \De^\prime=\fraco n,\De, ,\ 
\de^\prime=0,\ 0<\fraco \de,\De,\le \fraco 1,2,
 \big\}. 
\end{align*}
Кроме того, члены с $\fraco \de,\De,= \fraco 1,2,$ берутся в сумме с коэффициентом $\fraco 1,2,$.
Эта формула, в частности, даёт новое доказательство сравнений Рамануджана
$$
\tau (n)\equiv\sigma _{11}(n)(\bmod {691}).
$$

С помощью этих формул можно также найти  $\tau(6911) = -615012709514736031488$, причём оказывается, что
$\tau(6911)\equiv 1+6911^{11} (\bmod 691)$, но $\tau(6911)\not\equiv 1+6911^{11} (\bmod 691^2)$.

\subsubsection*{Вычисление с PARI-GP}
(см. также 
\begin{verbatim}
http://www.research.att.com/~njas/sequences/A000594, и
D. H. Lehmer, Tables of Ramanujan's function tau(n), Math. Comp., 24 (1970), 495-496.)

Программа на PARI-GP: 
{
m=11;n=2 ;si(n,m)= p=0; fordiv(n,d, p+= d^m); p \\ сумма степеней делителей
}
{
s1(n)=d3=1; vd1=[]; c1=0;   
\\ первая часть суммы (с \Delta>\delta>0,\ \Delta'>\delta'>0 
for(d1=1,n-1, for(d2=1,n-1,if(n-d1*d2>0, 
fordiv(n-d1*d2, d3,if(((d3<d1)& ((n-d1*d2)/d3<d2)),
d4=(n-d1*d2)/d3; c1=c1+1;
vd1=concat(vd1,[[d1,d2,d3,d4,c1]]);  
print("Delta="d1,"\t", "Deltap="d2,"\t","delta="d3,"\t", " deltap=" (n-d1*d2)/d3,"\t",c1);
)))));vd1
}
{
s2(n)= c2=c1; vd2=[];fordiv(n, d1, d2=n/d1;d4=0;for(d3=0,d1/2, \\ вторая часть суммы (с \delta'=0)
if(d3==d1/2, c=1/2, c=1); c2=c2+1; 
vd2=concat(vd2, [[d1,d2,d3,d4,c, c2]]);
print("Delta="d1,"\t", "Deltap="d2,"\t","delta="d3,"\t", "deltap="d4,"\t",  c ,"\t",c2)
)) ; vd2
}	
{
tau(n)=s1(n); s2(n); lvd1=length(vd1); lvd2=length(vd2); sn=0;
for(i1=1,lvd1, sn+= 
	vd1[i1] [1]^2* vd1[i1] [3]^2*(vd1[i1] [1]^2- vd1[i1] [3]^2)^3);
for(i2=1,lvd2,sn+=   
	(vd2[i2] [5])*vd2[i2] [1]^2* vd2[i2] [3]^2*(vd2[i2] [1]^2- vd2[i2] [3]^2)^3);
si(n,11)-(691/18)*sn
}

gp > tau(100)
Delta=2 Deltap=34       delta=1  deltap=32      1  
Delta=2 Deltap=35       delta=1  deltap=30      2
Delta=2 Deltap=36       delta=1  deltap=28      3
Delta=2 Deltap=37       delta=1  deltap=26      4

...............................................................................

Delta=100       Deltap=1        delta=50        deltap=0        1/2     291

%3 = 37534859200 \\ \\ результат: tau(100)
gp > ##
  ***   last result computed in 160 ms.
\end{verbatim}
Много других методов см. в \cite{Sloane}, \cite{Leh70}.
Отметим открытую проблему (проблема Лемера) о том, что $\tau(n)$ не обращается в нуль.

Другая нтересная открытая проблема состоит  в построении
 полиномиального алгоритма вычисления 
$\tau (p )$  для простого числа  $p $.
Аналогичный результат известен для коэффициентов $a (p)$ производяшего ряда $f_E(z)$ эллиптической кривой $E$ над $\Q$ (``алгоритм Схоофа''). По теореме Уайлса, такой ряд является модулярной формой веса 2 относительно некоторой конгруэнц-подгруппы модулярой группы,
в то время,как  $\tau (p )$ являются коэффициентами модулярной формы веса 12 относительно полной
 модулярой группы.   
\section{Классические модулярные формы}\label{ii43.2.}
 вводятся как некоторые  голоморфные функции на верхней комплексной полуплоскости
 $ {\HH} = \{ z\in {\CC}\ |\ \Im\ z>0\}$,
 которую можно также рассмативать как однородное пространство группы
  $G(\R) = \GL_2(\R)$:
\begin{align}\label{ii43.1}
 {{\HH}}
 =
 \GL_2(\R) / \Or (2)\cdot Z,
\end{align}
где $Z = \{ \mat x,0,0,x,\ | x\in \R^\times\}$
центр группы $G(\R)$ а $\Or (2)$ ортогональная группа.
При этом группа  $\GL_2^+(\R)$
матриц $\ga = \mat a_\ga, b_\ga, c_\ga, d_\ga,$ с положительным определителем действует на
 ${{\HH}}$ дробно-линейными преобразованиями; 
  на левых смежных классах  (\ref{ii43.1}) это действие переходит в естественное действие
 групповыми сдвигами. 

Пусть $\Ga$ -- подгруппа конечного индекса в модулярной группе
$\SL_2(\Z)$.
\begin{definition}
Голоморфная функция на верхней комплексной полуплоскости
 $ {\HH} = \{ z\in {\CC}\ |\ \Im\ z>0\}$ $f : {{\HH}} \to {\CC}$
 называется модулярной форой целого веса  $k$ относительно $\Ga$,
 если выполнены следуюшие условия a) и b):
\medskip
\begin{description}
\item{a)} {\it Условие автоморфности}
\begin{align}\label{ii43.2}
 f((a_\ga z + b_\ga)/(c_\ga z + d_\ga))
 =
 (c_\ga z + d_\ga)^k f(z)
\end{align}
для всех $\ga \in \Ga$;

\item{b)} {\it Регулярность в вершинах}:
$f$ регулярна в вершинах
 $z \in \Q \cup i\infty$;
 это означает, что для каждого элемент $\si = \mat a, b, c, d, \in 
\SL_2(\Z)$
функция $(cz + d)^{-k}f\left(\frac {az + b} {cz + d} \right )$
 разлагается в ряд Фурье по неотрицательным степеням $q^{1/N} = e(z/N)$
для некоторого натурального числа  $N$.
Модулярная форма
$$
f(z) = \sum _{n=0}^\infty a(n) e(nz/N)
$$
 называется  параболической, если   $f$ 
 обращается в нуль во  всех вершинах

 (т.е. их разложения Фурье содержат лишь строго положительные степени $q^{1/N}$),
\end{description}
\end{definition}
 see \cite{Serre J.--P. (1970)},
\cite{Ma-Pa05}, глава 6.
Комплексное векторное пространство всех модулярных форм  (соотв. параболических) форм
веса $k$ относительно $\Ga$
обозначается $\Mc_k(\Ga)$ (соотв.  $\Sc_k(\Ga)$).

Фундаментальный результат теории модулярных форм
 утверждает, что эти пространства конечномерны.
Кроме того, имеем $\Mc_k(\Ga)\Mc_l(\Ga) \subset \Mc_{k + l}(\Ga)$.
Прямая сумма
 $$
 \Mc (\Ga)
 =
 \bigoplus_{k=0}^\infty \Mc _k(\Ga)
 $$
 является градуированной алгеброй над
 ${\CC}$
 с конечным числом образующих.

Пример модулярных форм относительно $\SL_2(\Z)$ веса $k 
\ge 4$
даётся  {\it рядами Эйзенштейна}
\begin{align}\label{ii43.3}
 G_k(z)
 =
 \sum_{m_1, m_2 \in \Z}\!\!\!\!{}^{\displaystyle {}'} 
(m_1+ m_2 z)^{-k}
\end{align}
 (прим означает, что 
 $(m_1, m_2) \not= (0, 0)$).
Для этих рядов условие автоморфности  (\ref{ii43.2})
 непосредственно выводится из определения.
Имеем $G_k(z) \equiv 0$ для нечётных  $k$ и 
\begin{align}\label{ii43.4}
 G_k(z)
 =
 \frac{ 2(2\pi i)^k}{ (k-1)!}
 \left[ - \frac {B_k} {2k}  + \sum_{n=1}^\infty \si_{k-1}(n) e(nz) 
\right ],
\end{align}
где $\si_{k-1} (n) = \sum_{d|n}d^{k-1}$
и $B_k$ обозначает  $k^{ \hbox{е}}$
 {\it число Бернулли}.

Градуированная алгебра $\Mc (\SL_2(\Z))$    
 изоморфна кольцу многочленов от независимых переменных  $G_4$ и $G_6$.

\subsection{Фундаментальная область модулярной группы}\label{ens1.6.1}
Пусть
 $S=\mat 0, -1, 1, 0, $ et $T = \mat 1, 1, 0, 1, $. Имеем 
$$
S(z)=-z^{-1},\ T(z)=z+1.
$$
С другой стороны, пусть
 $D$ подмножество ${\HH}$ состояшее из точек  $z$
таких, что $|z|>\ge 1$ и $|\Re (z)|\le 1/2$.
Мы увидим, что 
 $D$ является  {фундаментальная областью } для действия модулярной группы 
 $\Ga(1)=\SL(2,\Z)$ на ${\HH}$, т.e. естественное отображение проекции 
$ D\to \Ga(1)\bs {\HH} $ сюръективно, а его ограничение на внутренность $D$ инъективно.
В то же время мы видели, что $S$ и $T$ порождают  $\Ga(1) = \SL(2,\Z)$. 

\begin{theorem}\label{fund}
1) Для всех  $z\in {\HH}$ существует матрица
 $\ga \in \Ga(1)$, такая, что  $\ga (z)\in D$. 

\noindent 2) Предположим, что две различные точки  $z, z^\prime \in D$
 эквивалентны при действии   $\Ga(1)$. 
Тогда или  $\Re(z)=\pm 1/2$ и $z=z^\prime+1$,
или $|z|=1$ и $z^\prime = -1/z$.

\noindent 3) Пусть $z\in D$, и пусть $St(z)=\{ \ga\in \Ga(1)\ |\ 
\ga(z)=z\}$ стабилизатор точки $z$ в $\Ga(1)$. 
Тогда имеем $St(z)=\{\pm 1\}$ за исключением трёх следующих случаев:

\noindent $z=i$, при этом $St(z)$ группа порядка  4 порождённая $S$;

\noindent $z=\rho =e^{2\pi i/3}$, при этом $St(z)$ группа порядка
 6 порождённая элементом $ST=\mat 0, -1, 1, 1, $;

\noindent $z=-\overline  \rho =e^{\pi i/3}$, при этом $St(z)$ группа порядка
 6 порождённая элементом $TS=\mat 1, -1, 1, 0, .$
\end{theorem}
\begin{figure}[h]
\begin{center}
\includegraphics[width=10cm,height=10cm]{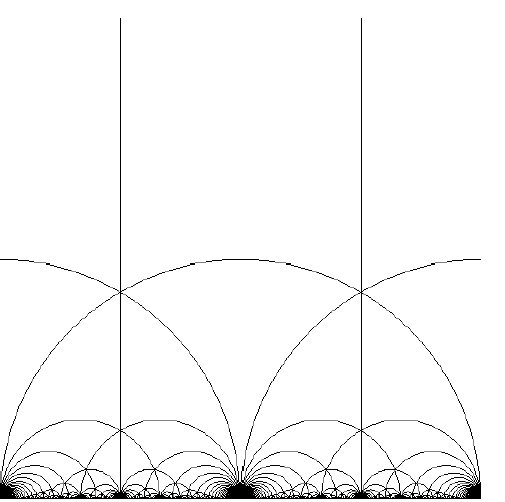}
\caption{\label{sl2z}Действие группы  $\SL(2,\Z)$.} 
На рисунке  \ref{sl2z} 
 представлено действие группы $\SL(2,\Z)$ на верхней комплексной полуплоскости. 
\end{center}
\end{figure}


Множество  ${{\HH}}/\SL_2(\Z)$ можно отождествить с можеством классов изоморфизма 
эллиптических кривых над ${\CC}$:
точке $z \in {{\HH}}$ сопоставляется комплексный тор 
 ${\CC}/(\Z + z\Z)$ который аналитически изоморфен римановой поверхности эллиптической кривой,
 записанной в форме Вейерштрасса:
\begin{align}\label{ii43.5}
 y^2
 =
 4x^3 - g_2(z)x - g_3(z)
\end{align}
где $g_2 = 60 G_4(z)$, $g_3(z) = 140G_6(z)$.

При замене $z$ на $\ga (z)$
для $\ga = \mat a_\ga, b_\ga, c_\ga, d_\ga ,\in  \SL_2(\Z)$
 решётка $\La_z = \Z + z\Z$ заменится на 
$$
 \La_{\ga(z)}
 =
 \Z + \ga (z)\Z
 =
 (cz + d)^{-1}(\Z + z\Z)
 =
 (cz + d)^{-1}\La_z, 
$$
 а кривая  (\ref{ii43.5}) примет каноническую форму Вейерштрасса с коэффициентами
$$
 g_2(\ga (z)) = (cz + d)^4g_2(z), \ \ \
 g_3(\ga (z)) = (cz + d)^6g_3(z).
$$
Дискриминант кубического многочлена справа 
(\ref{ii43.5})
 является модулярной параболической формой веса 12 относительно группы $\Ga = \SL_2(\Z)$:
\begin{align}\label{ii43.6}
 & 2^{-4}(g_2^3 - 27g_3^2)
 = \\ \nonumber &
 2^{-4}(2\pi)^{12}q \prod_{m=1}^\infty (1-q^m)^{24}
 = 
 2^{-4}(2\pi)^{12} \sum_{n=1}^\infty  \tau (n) q^n,
\end{align}
где $\tau (n)$ функция Рамануджана. 
При этом функция
\begin{align}\label{ii43.7}
 j(z)
 =
 1728 \frac {g_2^3}{ g_2^3 - 27g_3^2}
 =
 \frac 1 q  + 744 + \sum_{n=1}^\infty  c (n) q^n
\end{align}
мероморфна на  ${{\HH}}$ и в $\infty$,
 и не меняется при дробно-линейных преобразованиях с матрицами из  $\Ga = \SL_2(\Z)$.
 Эта функция доставляет важный пример {\it модулярной функции} и называется  {\it модулярным инвариантом}


\subsection{Модулярные формы как вычислительное средство решения задач арифметики}
Таким образом, мы можем рассматривать 
модулярные формы как \\
1) {\it степенные ряды} 
$\ds f=\sum_{n=0}^\infty a_n q^n\in {\mathbb C}[[ q]]$ и как \\
2) {\it голоморфные функции} 
{\it на верхней полуплоскости} 
$$ \HH = \{ z\in {\mathbb C}\ |\ \Im\ z>0\},$$ 
где $q=\exp(2\pi i z)$, 
$z\in \HH$,  и рассмотрим 
$L$-функцию 
$\ds L(f, s,\chi)=\sum_{n=1}^\infty \chi(n)a_n n^{-s}$ 
 для любого характера Дирихле 
$\chi:(\Z/N\Z)^*\to {\mathbb C}^*$, например, для символа Якоби  $\chi(n)=\left({n\over N}\right)$.

\subsection*{Ещё один метод вычисления функции Рамануджана:}
Положим
$ \ds
h_k:=\sum_{n=1}^\infty \sum_{d|n}d^{k-1}q^n
=\sum_{d=1}^\infty \frac {d^{k-1}q^d}{1-q^d}.
$ 

Доказывается: 
$\De=(E_4^3-E_6^2)/1728$, где
$E_4= 1+240h_4$ и
$ E_6= 1 - 504h_6$: 

\medskip
\noindent
{\bf Вычисление с PARI-GP}

\noindent
(см. \cite{BBBCO}).

\noindent
$ \ds
h_k:=\sum_{n=1}^\infty \sum_{d|n}d^{k-1}q^n
=\sum_{d=1}^\infty \frac {d^{k-1}q^d}{1-q^d}
\Longrightarrow
$ 

\begin{verbatim}
gp > h6=sum(d=1,20,d^5*q^d/(1-q^d)+O(q^20))
gp > h4=sum(d=1,20,d^3*q^d/(1-q^d)+O(q^20)
gp > Delta=((1+240*h4)^3-(1-504*h6)^2)/1728
\end{verbatim}

\begin{verbatim}
q - 24*q^2 + 252*q^3 - 1472*q^4 + 4830*q^5 - 6048*q^6 - 16744*q^7 
+ 84480*q^8 - 113643*q^9 - 115920*q^10 + 534612*q^11 
- 370944*q^12 - 577738*q^13 + 401856*q^14 + 1217160*q^15 
+ 987136*q^16 - 6905934*q^17+ 2727432*q^18 + 10661420*q^19 + O(q^20)
\end{verbatim}

\noindent
{\bf Сравнение Рамануджана: }
$\ds \tau (n) \equiv \sum_{d|n}d^{11} \ \bmod \ 691:$

\begin{verbatim}
gp > (Delta-h12)/691
%10 = -3*q^2 - 256*q^3 - 6075*q^4 - 70656*q^5 - 525300*q^6
 - 2861568*q^7 - 12437115*q^8 - 45414400*q^9
 - 144788634*q^10 - 412896000*q^11 - 1075797268*q^12
 - 2593575936*q^13 - 5863302600*q^14 - 12517805568*q^15
 - 25471460475*q^16 - 49597544448*q^17
 - 93053764671*q^18 - 168582124800*q^19 + O(q^20)
\end{verbatim}
\subsection*{Вот ешё три программы вычисления $\tau(n)$ (см. \cite{Sloane})}
\begin{verbatim}
 PROGRAM   	

(MAGMA) M12:=ModularForms(Gamma0(1), 12); t1:=Basis(M12)[2]; 
PowerSeries(t1[1], 100); Coefficients($1);

(PARI) a(n)=if(n<1, 0, polcoeff(x*eta(x+x*O(x^n))^24, n))

(PARI) {tau(n)=if(n<1, 0, polcoeff(x*(sum(i=1, (sqrtint(8*n-7)+1)\2,
(-1)^i*(2*i-1)*x^((i^2-i)/2), O(x^n)))^8, n));} 
gp > tau(6911)
%3 = -615012709514736031488
gp > ##
  ***   last result computed in 3,735 ms.
\end{verbatim}

\subsection*{Схема применения  модулярных форм для решения задач теории чисел:}
\hskip-1cm
\begin{tabular}{cccc} 
\begin{tabular}{l}
Производящая\\ 
функция \\
$f=\sum_{n=0}^\infty a_n q^n$ \\
$\in {\mathbb C}[[ q]]$\\
для арифметической \\
функции $n\mapsto a_n$, \\
например $a_n=p(n)$
\end{tabular}
&
$\rightsquigarrow$
\begin{tabular}{l}
Выражение через\\
модулярную форму,\\
например \\
$\ds\sum_{n=0}^\infty p(n) q^n$\\
$=(\Delta/q)^{-1/24}$
\end{tabular}
&
$\rightsquigarrow$ 
\begin{tabular}{l}
Число\\
(ответ)
\end{tabular}
\\ 
\begin{tabular}{l}
{\bf Пример 1} (см. \rm \cite{Chand70}): \\ 
(Харди-Рамануджан)
\end{tabular} 
&
$\uparrow$
&
$\uparrow$
\\
\tiny
\begin{tabular}{l}
$\ds p(n)=\frac{e^{\pi \sqrt{2/3({n-1/24})}}} {4\sqrt{3}\lambda_n^2}$ \\  
$+O(e^{\pi \sqrt{2/3}\lambda_n}/ \lambda_n^3),$\\
$\lambda_n=\sqrt{n-1/24}$,
\end{tabular}
&
\begin{tabular}{l}
Хорошие базисы \\
конечномерность\\
много соотношений \\ 
и тождеств
\end{tabular}
&
\begin{tabular}{l}
Значения \\
$L$-функций, \\ 
сравнения, \\
\dots
\end{tabular}
\end{tabular}

\noindent
{\bf Пример 2}
(см. в \cite{Ma-Pa05}, главы 6 и 7): теорема Ферма-Уайлса, гипотеза Бёрча-Суиннертона-Дайера, \dots


 \section{Ряды Эйзенштейна и сравнения для функции Рамануджана.
 }
\subsubsection*{Ряды Эйзенштейна и их разложение Фурье.}\label{ens1.5}
Пусть $k>2$. Для решётки  $\La\subset {\CC}$ положим
$$
G_k(\La)=\sum_{l\in \La}l^{-k}=\sum_{m,n}\!{}^{\displaystyle {}^\prime}(m\w_1+n\w_2)^{-k}, \ \ 
\La=\langle \w_1, \w_2\rangle,
$$
Этот ряд сходится абсолютно для  $k> 2$.

\begin{proposition}  \ 
\noindent
{\hbox{(а)}}  Имеем 
$$
G_k(z)= \sum_{m,n \in \Z}\!\!\!\!{}^{\displaystyle {}'} 
(mz + n)^{-k}\in  \Mc_k(\Ga(1));
$$

\medskip\noindent
{\hbox {(б)}}
$$
G_k(z)=
2\z (k)\left[ 1-\fraco 2 k, B_k, \sum_{n=1}^\infty \si_{k-1}(n)q^n\right]
=:
2\z (k)E_k(z),
$$
где $q=e(z)=\exp (2\pi i z)$, $B_k$ числа Бернулли определённые разложением 
$$
\fraco x, e^x-1, =\sum_{k=0}^\infty B_k\fraco x^k, k!, 
$$
 \end{proposition} 
Вот несколько численных значений: 
\begin{align*}
&
B_0=1, \ B_1=-\fraco 1, 2,,  \ B_2=\fraco 1, 6, ,  \ B_3=B_5=\cdots=0, \ 
B_4=-\fraco 1, 30,, \ B_6= \fraco 1, 42, , \cr &\ B_8=-\fraco 5, 66, ,\  
B_{12}=\fraco 691, 2730,,\ B_{14} =-\fraco 7, 6, ,\ B_{16}=\fraco 3617, 510,, \ 
B_{18} = -\fraco 43867, 798, , \  \dots.
\end{align*}
Имеем $\displaystyle{\z(k)=-\fraco (2\pi i)^k, 2, \fraco B_k, k!,} $, 
$$
G_k(z)= \fraco (2\pi i)^k, (k-1)!, 
\left[ -\fraco B_k, 2 k, + \sum_{n=1}^\infty
\si_{k-1}(n)q^n\right]=:\fraco (2\pi i)^k, (k-1)!, \GG_k(z).
$$

\medskip
\noindent
{\bf Примеры.}
\begin{align*}
&
E_4(z)=1+240\sum_{n=1}^\infty\si_{3}(n)q^n \in \Mc_4(\SL(2, \Z)),\cr &
E_6(z)=1-504\sum_{n=1}^\infty\si_{5}(n)q^n \in \Mc_6(\SL(2, \Z)),\cr &
E_8(z)=1+480\sum_{n=1}^\infty\si_{7}(n)q^n\in \Mc_8(\SL(2, \Z)),\cr &
E_{10}(z)=1-264\sum_{n=1}^\infty\si_{9}(n)q^n\in \Mc_{10}(\SL(2, \Z)),\cr &
E_{12}(z)=1+\fraco 65520, 691, \sum_{n=1}^\infty\si_{11}(n)q^n
\in \Mc_{12}(\SL(2, \Z)),\cr &
E_{14}(z)=1 - 24\sum_{n=1}^\infty\si_{13}(n)q^n\in \Mc_{14}(\SL(2, \Z)).
\end{align*}
\medskip
\noindent
{\it Доказательство. }Автоморфность ясна, поскольку
 $G_{k}(\lambda \La)=\lambda^{-k}G_k(\La)$ поэтому $G_k$ 
является {\it однородной функцией решётки} степени однородности  $-k$, и
\begin{align*}&
G_k(z)= G_{k}(\La_z),\  G_k(\gamma z)= G_{k}(\La_{\gamma z})=
G_{k}(\langle 1, {\frac {az+b}{cz+d}\rangle})
\\ &
=G_{k}(({cz+d})^{-1}\langle {cz+d}, { {az+b}\rangle})=
({cz+d})^{k}G_{k}(\langle {cz+d}, { {az+b}\rangle})=({cz+d})^{k}G_{k}(\La_z)=
({cz+d})^{k}G_{k}(z),
\end{align*}
поскольку $\langle {cz+d},  {az+b}\rangle=\langle 1,  z\rangle $ для всех $\mat a,b,c,d,\in SL_2(\Z)$.

Для нахождения разложения Фурье используется известное разложение 
синуса в бесконечное произведение:
\begin{align}\label{ens5.1}
\sin (\pi a)=\pi a \prod_{n=1}^\infty\left (1-\fraco a^2, n^2, \right).
\end{align}
Логарифмическая производная 
(\ref{ens5.1}) даёт
\begin{align}\label{ens5.2}
\pi {\rm ctg} \pi a =\fraco 1, a, +\sum_{n=1}^\infty\left 
(\fraco 1, a+n, -\fraco 1,
a-n, \right ).
\end{align}
Заметим, что
\begin{align}\label{ens5.3}
\pi i \fraco e^{\pi i a }+e^{-\pi i a },e^{\pi i a }-e^{-\pi i a }, =
\pi i + \fraco 2 \pi i , e^{2\pi i a }-1, = 
\pi i - 2 \pi i \sum_{n=1}^\infty e^{2\pi i na },
\end{align}
и положим $x=2 \pi i a$; отсюда
$$
\fraco x, 2, +\fraco x, e^x-1, =1+\sum_{n=1}^\infty \fraco 2x^2, x^2-(2\pi i n)^2,, 
$$
где
\begin{align*}
&
\sum_{k=0}^\infty B_k\fraco x^k, k!,+\fraco x, 2, =1-\sum_{n=1}^\infty \fraco
2\left( \displaystyle{\fraco x, 2 \pi i n,} \right)^2, -\left(
\displaystyle{\fraco x, 2 \pi i n, }\right)^2+1, =
\cr & 1-2\sum_{n=1}^\infty \sum_{k=2k^\prime\ge 2}
\left( \displaystyle {\fraco x, 2 \pi i n,} \right)^k=
1-2\sum_{k=2k^\prime\ge 2}\displaystyle {\fraco \z(k), (2\pi i)^k,} x^k.
\end{align*}
Это непосредственно даёт
\begin{align}\label{ens5.4}
\z (k)=-\fraco (2\pi i)^k, 2, \fraco B_k, k!, ,
\end{align}
в частности, 
$$
\z (2)=\fraco \pi ^2, 6, , \z(4)=\fraco \pi ^4, 90, .
$$
Чтобы доказать (б),  проводится 
дифференцирование
обеих частей  (\ref{ens5.3})
по переменной  $a$ $(k-1)$ раз:
\begin{align}\label{ens5.5}
-(2\pi i)^k \sum_{n=1}^\infty n^{k-1}e^{2\pi i n a}=
(-1)^{k-1}(k-1)! \sum_{n\in \Z}(a+n)^{-k}, \ \ (k\in 2\Z, k\ge 2).
\end{align} 
Положим $a=mz$, тогда
\begin{align}\label{ens5.6}
\fraco (2\pi i)^k, (k-1)!, \sum_{n=1}^\infty n^{k-1}e^{2\pi inmz}=
\sum_{n\in \Z}(mz+n)^{-k}. 
\end{align}
Если теперь $k>2$, то можно просуммировать по  $m$
от 1 до $\infty$.
В результате этого получим
\begin{align}\label{ens5.7}
G_k(z)=
2\z (k) + 2\sum_{m=1}\sum_{n=-\infty}^{\infty}(mz+n)^{-k}=
2\z (k)\left[ 1-\fraco 2 k, B_k, \sum_{m, d=1}^\infty d^{k-1}q^{md}\right].
\end{align}
Отметим, что двойной ряд в  (\ref{ens5.7}) абсолютно сходится при
 $k>2$
 но ряд   (\ref{ens5.7}) имеет смысл и при  $k=2$ 
 как условно сходящийся ряд. 
Доказательство завершается подстановкой  (\ref{ens5.4}) в (\ref{ens5.7}).

\

\begin{theorem}\label{ens1.5.2}  
Пусть
 $\displaystyle{\De(z)=q\prod_{m\ge 1}}(1-q^m)^{24}$. 
Тогда имеем
$$
\De (-z^{-1})=z^{12}\De(z).
$$
\end{theorem} 
(см. также \cite{Serre J.--P. (1970)}). 

\medskip
\noindent
{\it Доказательство.} Положим
$$
E_2(z)=1-24\sum_{n=1}^\infty\si_{1}(n)q^n.
$$
Имеем
\begin{align*}
&
\fraco d, dz, \log (\De(z))=\fraco d, dz, \log q+
24\sum_{m=1}^\infty \fraco d, dz, \log(1-q^m)=
\cr &
2 \pi i (1-24 \sum_{m=1}^\infty mq(1-q^m)^{-1})=
2\pi i E_2(z),\ \ \fraco dq, dz, =2\pi i q.
\end{align*}
Достаточно доказать следующее предложение:

\begin{proposition}\label{ens1.5.3}
$$
z^{-2}E_2(-z^{-1})=E_2(z)+\fraco 12, 2\pi i z, .\eqno(5.8) 
$$
\end{proposition}

{\it Доказательство предложения.} Используется ряд (\ref{ens5.7}) с $k=2$
сходящийся условно:
\begin{align*}
&
E_2(z)=
\fraco 1,
2\z(2),\sum_{m=-\infty}^\infty\left (\sum_{n=-\infty\atop n\not = 0}^{\infty}
(mz+n)^{-2}\right )=
\cr &1+\fraco 3, \pi^2, \sum_{m\not = 0}\left 
(\sum_{n=-\infty}^{\infty}
(mz+n)^{-2}\right )=
1+\fraco 6, \pi^2, \sum_{m=1}^\infty\left 
(\sum_{n=-\infty}^{\infty}
(mz+n)^{-2}\right ).
\end{align*}
Для фиксированного  $m$ имеем
$$
\sum_{n=-\infty}^{\infty}
(mz+n)^{-2}=1- \fraco 4, B_2, \sum_{d=1}^\infty dq^{md}=
1-24\sum_{n=1}^\infty\si_{1}(n)q^n.
$$
Выполним подстановку
$$
z^{-2}E_2(-z^{-1})=\fraco 1,
2\z(2),\sum_{m=-\infty}^\infty\left (\sum_{n=-\infty\atop n\not = 0}^{\infty}
(-m+nz)^2\right)
=
1+\fraco 3, \pi^2, \sum_{n=-\infty}^{\infty}\sum_{m\not =0}
(mz+n)^{-2}.
$$
Если положить $a_{m, n}=(mz+n)^{-2}$, то доказательство сводится к проверке равенства
$$
-\sum_m\sum_n a_{m, n} + \sum_n\sum_m a_{m, n}=\fraco 12, 2\pi i z, .
$$
Для его доказательства вводится поправочный член
\begin{align}\label{ens5.9}
b_{m, n}(z)=\fraco 1, (mz+n-1)(mz+n), =\fraco 1, (mz+n-1), -\fraco 1, (mz+n),
\end{align}
Получается модифицированный ряд 
\begin{align}\label{ens5.10}
\tilde E_2(z)=1+\fraco 3, \pi^2,  \sum_{m\not = 0}
\sum_{n=-\infty}^{\infty}\left((mz+n)^{-2}-b_{m, n}(z)\right)
\end{align}
который уже абсолютно сходится
поскольку
$$
(mz+n)^{-2}-((mz+n-1)(mz+n))^{-1}=
(mz+n)^{-2}(mz+n-1)^{-1}.
$$
С другой стороны,
\begin{align*}
&
\tilde E_2(z)=\cr &
1+\fraco 3, \pi^2,  \sum_{m\not = 0}\left 
(\sum_{n=-\infty}^{\infty}(mz+n)^{-2}\right )+
\fraco 3, \pi^2,  \sum_{m\not = 0}\sum_{n=-\infty}^{\infty}
\left (\fraco 1, (mz+n), -\fraco 1, (mz+n-1), \right),
\end{align*}
и последняя сумма преобразуется в нуль, поэтому
$$
\tilde E_2(z)=E_2(z).
$$
Изменение порядка суммирования в (\ref{ens5.10}) обосновано в силу абсолютной сходимости,
откуда
\begin{align*}
&
\tilde E_2(z)=1+\fraco 3, \pi^2, 
\sum_{n=-\infty}^\infty
\sum_{m\not =0}\big((mz+n)^{-2}-b_{m, n}(z)\big)=\cr &
z^{-2}E_2(-z^{-1})-\fraco 3, \pi^2, 
\sum_{n=-\infty}^\infty
\left(\sum_{m\not =0}b_{m, n}\right).
\end{align*}
Остаётся вычислить последнюю сумму:
$$
\sum_{n=-\infty}^\infty
\left(\sum_{m\not =0}b_{m, n}\right)=
\lim_{N\to \infty}\sum_{n=-N+1}^{n=N}\left(\sum_{m\not =0}b_{m, n}\right).
$$
Однако
$$
\sum_{m\not =0}(mz-n)^{-2}=\fraco 1, z^2, \sum_{m\not =0}
(n/z-m)^{-2}=-\fraco 1, n^2, -\fraco 4\pi^2, z^2, 
\sum_{d=1}^\infty de^{-2\pi i nd(1/z)}
$$
поэтому для всех $z$ 
внешняя сумма сходится абсолютно, и преобразуется в 
\begin{align*}
&\sum_{m\not =0}\left(\sum_{n=-N+1}^{n=N}b_{m, n}\right)
=  \sum_{m\not =0}\left (\fraco 1, (mz-N), -\fraco 1, (mz+N), \right)=
\cr &
\fraco 2, z, \sum_{m=1}^\infty\left (\fraco 1, (-N/z+m), +
\fraco 1, (-N/z-m), \right)=\fraco 2, z,\left (\pi {\rm ctg} \left(-\fraco \pi N, z,
\right)+\fraco z, N, \right) \to -\fraco 2\pi i, z, 
\end{align*}
при $N\to \infty$, $z\in {\HH}$,  откуда следует и предложение  \ref{ens1.5.3}, и теорема
 \ref{ens1.5.2}.

\subsection{Структура пространств модулярных форм относительно $\SL_2(\Z)$.}\label{1.6ens}
 (см. также \cite{Serre J.--P. (1970)}, pp.127--178).

Пусть $f$ ненулевая мероморфная функция на  ${{\HH}}$, и пусть  $p$ некоторая точка в ${\HH}$.
Назовём порядком  $f$ в $p$, и обозначим его через  $v_p(f)$,
целое число  $n$ такое, что функцмя $f/(z-p)^n$ голоморфна и необращается в нуль в точке $p$.

Пусть $f$ модулярная функция веса $k$, то равенство  
$$
f(z)=(cz+d)^{-k}f\left(\fraco az+b, cz+d, \right)
$$
показывает, что  $v_p(f)=v_{\ga(p)}(f)$ для всех $\ga\in \Ga=\Ga(1)$; 
другими словами, $v_p(f)$ зависит только от образа  $p$ в $\Ga\bs H$.
Больше того, можно определить и  $v_\infty(f)$ как порядок относительно  $q=0$ функции 
$\tilde f(q)=f(z)$ ассоциированной с  $f$.
Положим $e_p=2$ (соотв.$e_p=3$) если $p$ эквивалентна относительно  $\Ga$ точке $i$ (соотв. точке
$\rho$), и $e_p=1$ в противном случае.

\begin{proposition}[о степени дивизора модулярной формы 
 ]\label{ens1.6.4} 
Пусть $f$ не\-ну\-ле\-вая модулярная функция веса $k$ относительно $\Ga(1)$.
Имеем
$$
v_\infty(f)+\sum_{p\in \Ga(1)\bs {\HH}}\fraco 1, e_p, v_p(f) = \fraco k, 12, 
$$
\end{proposition}

\medskip
\noindent
[Можно также записать этот результат в виде:
$$
v_\infty(f)+\fraco 1, 2, v_i(f) +\fraco 1, 3, v_\rho(f)
+\sum_{p\in \Ga(1)\bs {\HH}}\!\!\!\!\!{}^{{{}^{*\prime}}} v_p(f)= \fraco k, 12,,
$$
где символ $\displaystyle{\sum_{p\in \Ga(1)\bs {\HH}}\!\!\!\!\!{}^{{{}^{*\prime}}}}$ 
означает суммирование по всем классам точек  $\Ga(1)\bs {\HH}$, отличным от классов точек
$i$ et de $\rho$].

\noindent
Естественное доказательство этого факта использует структуру римановой поверхности на
 $\Ga(1)\bs \overline  {\HH}$, 
 где  $\overline  {\HH}={\HH} \cup \Q \cup \infty$. 

\begin{theorem}[о функции Рамануджана $\De$ и рядах Эйзенштейна]\label{ens1.6.5} 
(i)  Имеем $\Mc_k(\Ga(1))=0$ pour $k<0$ и $k=2$.

\medskip
\noindent
(ii)  Для $k=0, 4, 6, 8, 10$  пространство $\Mc_k(\Ga(1))$  имеет размерность  1
с базисом $1, E_4, E_6, E_8, E_{10}$;
при этом $\Sc_k(\Ga(1))=0$.

\medskip
\noindent
(iii) Умножение на $\De$ определяет изоморфизм 
$\Mc_{k-12}(\Ga(1))$ на $\Sc_k(\Ga(1))$.
\end{theorem}

\begin{theorem}[размерности пространств модулярных форм
для $\SL(2,\Z)$]\label{ens1.6.6}  
\ 

\hbox{(а)} 
$$
\dim \Mc_k(\Ga(1))=\begin{cases}
\left[\fraco k, 12, \right ],
&k \equiv 2 (\bmod 12), k\ge 0,\cr
0, &k\equiv 1 (\bmod 2),\cr
\left[\fraco k, 12, \right ]+1,
&k \not \equiv 2 (\bmod 12), k\ge 0, k\in 2\Z.
\end{cases}
$$
$$
\dim \Sc_k(\Ga(1))=
\begin{cases}
\left[\fraco k, 12, \right ]-1,
&k \equiv 2 (\bmod 12), k\ge 12,\cr
0, &k\equiv 1 (\bmod 2),\cr
\left[\fraco k, 12, \right ],
&k \not \equiv 2 (\bmod 12), k\ge 0, k\in 2\Z.
\end{cases}
$$

\medskip\noindent
\hbox{(б)}  Произведения
$$
\{ E_4^\al E_6^\beta\ |\ 4\al+6\beta=k,\ \al, \beta\ge 0, \al, \beta\in\Z\}
$$
образуют базис пространства $\Mc_k(\Ga(1))$
\end{theorem}

\medskip
\noindent
{\it Доказательство } непосредственно следует из  \ref{ens1.6.5}.

\begin{corollary}
Справедливо равенство
$$
\De(z)=\displaystyle{\fraco 1, 1728, (E_4^3-E_6^3)}.
$$
\end{corollary}

\medskip
\noindent
Действительно, $\De(z)\in \Sc_{12}(\Ga(1))$, и в силу  \ref{ens1.2.6} имеем $\dim \Sc_{12}(\Ga(1))=1$, 
остаётся заметить, что функция  $\displaystyle{\fraco 1, 1728,
(E_4^3-E_6^3)}$ также принадлежит одномерному пространству $\Sc_{12}(\Ga(1))$, так как  обе функции $E_4^3, E_6^3$
имеют  коэффициент при $q$, равный 1.

\subsection{Приложение: доказательство сравнения Рамануджана}\label{ens1.6.7}
\begin{align}\label{Ram}
\tau (n)\equiv \si_{11}(n)\bmod 691.
\end{align}
Действительно, 
$$
E_6^2(z)-\left(1-504\sum_{n=1}^\infty \si_5(n)q^n\right)^2 \in \Z [[q]],
$$
поэтому можно разложить  $E_6^2(z)$ в базисе $\{E_{12}$, $\De\}$ пространства
$\Mc_{12}(\Ga(1))$ размерности 2:
$E_6^2=E_{12}+\al \De$, где 
$$
1-1008q + \dots = 1+ \fraco 65 520, 691, q +\dots + 
\al q + \dots, \ \ 
$$
и $\dots = {\oc} (q^2)$.
Поэтому
$$
\al=-1008-\fraco 65 520, 691,=\fraco a, 691, \ \ \hbox{где}\ a\equiv -65 520
(\bmod 691),
$$
и из разложения выводится. что 
$$
\fraco 65 520, 691,\si_{11} (n) + \fraco a, 691, \tau(n) \in \Z, \ \ \ \hbox{где}\
\  65 520(\si_{11}(n) -\tau(n)) \equiv 0 (\bmod 691),
$$
откуда вытекает сравнение (\ref{Ram}).

\section{Числа Бернулли и сравнения Куммера}
\subsection{Сравнения для коэффициентов рядов Эйзенштейна}
Приведем пример сравнений между коэффициентами модулярных форм по модулю $p^n$.

Для этого  рассмотрим ещё одну нормализацию рядов Эйзенштейна, заданную так,
что коэффициенты Фурье $a(n)$ задают ряд Дирихле с эйлеровским произведением, при этом $a(1)=1$: 
$$ 
{\GG}_k={\zeta(1-k)\over
2}E_k=-{B_k\over
2k}+\sum\limits_{n=1}^\infty
\sigma_{k-1}(n)q^n  = \sum_{n=0}^\infty a(n)q^n \Rightarrow \sum_{n=1}^\infty a(n)n^{-s}= \zeta(s)\zeta(s+1-k),
$$ 
а также $p$-нормализацию
$$ 
{\GG}_k^*(z)={\GG}_k(z)-p^{k-1}{\GG}_k(pz).
$$
Тогда
 
\begin{eqnarray} & & \nonumber \ds
{\GG}_k^*={\zeta^*(1-k)\over
2}+\sum\limits_{n=1}^\infty
\sigma_{k-1}^*(n)q^n,\,\sigma_{k-1}^*(n)=
\sum\limits_{{\scriptstyle d\vert n\atop
\scriptstyle (d,p)=1}}d^{k-1},\quad\hbox{где}\\ & & \nonumber \ds
\zeta^*(s) =\zeta(s)(1-p^{-s})=
\sum\limits_{{\scriptstyle n=1\atop
\scriptstyle (p,n)=1}}n^{-s}\quad
\hbox{\vtop{\hbox{обозначает дзета-функцию Римана}
\hbox{с удалённым эйлеровским  $p$-множителем.}}}
\end{eqnarray}
$$ {\GG}^*_k= \sum_{n=0}^\infty a^*_k(n)q^n \Rightarrow \sum_{n=1}^\infty a^*_k(n)n^{-s}= \zeta(s)\zeta^*(s+1-k).
$$ 

\normalsize
\begin{theorem}\label{THkor1.2}
\begin{description}
\item {а)}\enspace Пусть
$k\equiv k'\bmod{(p-1)}p^{N-1}$ тогда  ${\GG}_k^*\equiv {\GG}_{k'}^*\bmod{p^N}$ в
$\dbQ{[\![q]\!]}$ для всех $k\not\equiv 0 \bmod(p-1)$.

\item{б)}\enspace Пусть $k\equiv k'\bmod{(p-1)}p^{N-1}$, тогда для любых 
 $c\in\dbZ$, $(c,p)=1$,
$c>1$ имеем:  $(1-c^k){\GG}_k^*\equiv (1-c^{k'}){\GG}_{k'}^*
\bmod{p^N}$ (без ограничения на  $k$).

\item{в)}\enspace
Семейство классических модулярных форм 
$$
k\mapsto f_k=(1-c^k){\GG}_k^*
$$
 является $p$-адически непрерывным
 $\dbZ_p^*$ с параметром из множества
$$
\Pc =\{y\mapsto y^k, k\ge 4\}
$$ 
$p$-адических характеров группы $\dbZ_p^*$.
\end{description}
\end{theorem}

{\it Доказательство }{теоремы \ref{THkor1.2}}:
Утверждения   а) и б) следуют из в).
Для доказательства в) положим
$f_k=\sum\limits_{n\ge 0}a_k(n)q^n$ \\
\begin{description}
\item {\sf Случай $n>0$:}\enspace
Функции 
$$
k\mapsto a_k(n)=(1-c^k)\sum\limits_{{\scriptstyle
d\vert n\atop\scriptstyle (d,p)=1}}d^{k-1}
$$
являются $p$-адически непрерывными по их элементарному описанию (по сравнениям теоремы Эйлера);  

\item {\sf Случай $n=0$:}\enspace
$a_k(0)=(1-c^k)\zeta^*(1-k)$ рассматривается с помощью классических сравнений Куммера:
зафиксируем произвольное целое число $c\in\dbZ$, $(c,p)=1$,
$c>1$.
\end{description}

\begin{theorem}[Куммер]\label{kor1.3}
Пусть $\zeta_{(p)}^{(c)}(-k)=
(1-c^{k+1})(1-p^k)\zeta(-k)$, $k\ge 0$,
и пусть $h(x)=\sum\limits_{i}\alpha_i
x^i\in\dbZ[x]$ такой, что  $h(a)\equiv
0\bmod{p^N}$ для всех $\alpha\in\dbZ_p^*$.
Тогда $\sum\limits_{i}
\alpha_i\zeta_{(p)}^{(c)}(-i)
\equiv 0\bmod{p^N}$.
\end{theorem}


{\it Доказательство теоремы } \ref{kor1.3}
 использует суммы стененей
 $S_k(M)=\sum\limits_{n=1}^{M-1}n^k$, числа Бернулли $B_{k}$, и многочлены Бернулли
 $B_{k}(x)$:
$$
S_k(M)=\sum\limits_{n=1}^{M-1}n^k={1\over
k+1}[B_{k+1}(M)-B_{k+1}], \mbox{ где }  \sum\limits_{m=1}^\infty
\frac {B_k}{k!} t^k=\frac {t} {e^t-1}
\mbox{ и }B_{k}(x)=\sum_{i=0}^k{k \choose
 i}B_ix^{k-i}.
$$
Отсюда вытекает 
$$
B_{k} =
\lim_{N\rightarrow \infty} \frac 1 {p^N} S_{k}(p^N) \ 
$$ 
($p$--адический предел
явно вычсиляется по указанной формуле для  
$S_k(p^N)$ (в частности, $\frac 1 {p^N} S_{1}(p^N) = \frac {p^N (p^N-1)} {2p^N}\to
 -\frac 1 {2}=B_{1})$.

Далее, рассматривается регуляризованная сумма степеней 
$$
S_k^*(p^N)=\sum\limits_{{\scriptstyle
n=1\atop\scriptstyle p\nmid n}}^{p^N-1}n^k=
S_k(p^N)-p^kS_k(p^{N-1})\,\,, 
$$
 которая выражается через числа Бернулли в терминах $S_k(N)$
 по формуле  
$$
B_{k+1} =
\lim_{N\rightarrow \infty} \frac 1 {p^N} S_{k+1}(p^N). \ 
$$ 
Для всех  $n$ с $(p, n) = 1$ имеем сравнение
$h(n)  \equiv 0 (\bmod p^N)$, и
\begin{eqnarray*}
   & &\lim_{N\rightarrow \infty} \frac 1 {p^N} S_{k+1}^*(p^N)=
\lim_{{N}\rightarrow \infty} 
\frac 1  {p^N} [S_{k+1}(p^{N}) - p^{k+1}S_k(p^{{N}-1})] = \cr
  & &\lim_{{N}\rightarrow \infty} \frac 1  {p^N} S_{k+1}(p^{N}) - 
p^{k}\lim_{{N}\rightarrow \infty} \frac 1  {p^N}  S_{k+1}(p^{{N}-1}) =
(1-p^{k})B_{k+1}.\nonumber
\end{eqnarray*}

Подставим
$\ds \zeta(-k)=-{B_{k+1}\over k+1}$ тогда
\begin{align}\label{kor1.2}
\zeta_{(p)}^{(c)}(-k)=(c^{k+1}-1)(1-p^k){B_{k+1}
\over k+1}\equiv
{S_{k+1}(p^M)\over p^M}\cdot
{(c^{k+1}-1)\over k+1}\bmod{p^N}
\end{align}
(для достаточно большого  $M\ge N$).
Правая часть  (\ref{kor1.2}) преобразуется к виду 
\begin{align}\label{EQkor1.3}
\sum\limits_{{\scriptstyle n=1\atop
\scriptstyle p\nmid n}}^{p^M-1}{(cn)^{k+1}-n^{k+1}
\over p^M\cdot(k+1)}=
\sum\limits_{{\scriptstyle n=1\atop
\scriptstyle p\nmid
n}}^{p^M-1}{(cn)^{k+1}-n_c^{k+1}\over p^M\cdot
(k+1)}
\end{align}
где $n\mapsto n_c$ перестановка множества
$\{1,2,\ldots,p^M-1\}$ заданная $n_c\equiv
nc\pmod{p^M}$.
Подставим $cn=n_c+p^M t_n$, $t_n\in
\dbZ$ в (\ref{EQkor1.3}):
$$
{(nc)^{k+1}-n_c^{k+1}\over p^M\cdot (k+1)}
\equiv t_n\cdot n_c^k\bmod{p^M}
$$
поэтому
$ \ds
\zeta_{(p)}^{(c)}(-k)\equiv 
\sum\limits_{{\scriptstyle n=1\atop
\scriptstyle p\nmid n}}t\cdot 
n_c^k\bmod{p^M}
$
где $t_n=t(n,c)$ не зависит от  $k$.
Чтобы завершить доказательство, подставим это сравнение в линейную комбинацию из теоремы
\ref{kor1.3} 
используя {\sf предположение}
$$
h(x)\equiv 0\bmod p^N: \ \ 
\sum\limits_{i}\alpha_i\zeta_{(p)}^{(c)}
(-i)
\equiv 
\sum\limits_{{\scriptstyle n=1\atop
\scriptstyle p\nmid n}}t_n\cdot h(n_c)
\equiv 0\bmod{p^N}\,\,.\qed
$$

\begin{corollary}\label{cor1}
 ($p$-адическая непрерывность  $\zeta_{(p)}^{(c)}(-k)$ в прогрессиях по  $\bmod{(p-1})$).
Если  $h(x)=x^k-x^{k'}$, $k\equiv k'\bmod(p-1)
p^{N-1}$, то
$$
\zeta_{(p)}^{(c)}(-k)\equiv\zeta_{(p)}^{(c)}
(-k')\bmod{p^N}\,\,.
$$
\end{corollary}

{\it Доказательство  }{следствия \ref{cor1}}: По теореме Эйлера имеем 
$h(a)\equiv 0\bmod{p^N}$, поскольку
$a^{\varphi(p^N)}\equiv 1\pmod{p^N}$, $(a, p)=1$.

\subsection{$p$-адическое интегрирование и 
мера Мазура}

В $p$-адической теории интегрирования рассматривается {\it поле Тэйта},
$\CC_p = \widehat{\overline{\QQ}}_p$ (пополнение алгебраического замыкания поля  $p$-адических чисел        $\QQ_p)$, которое служит аналогом поля комп\-лекс\-ных чисел, так как $\CC_p$ алгебраически замкнуто,
 и является топологически полным мет\-ри\-чес\-ким пространством с нормой $|\cdot|_p$, $|p|_p=\frac 1 p$.

Пусть $R$  любое замкнутое подкольцо в $\CC_p$, $\Mc$ -- топологический $R$-модуль и $\Cr(Y,R)$ -- топологический модуль всех $R$-значных непрерывных функций на проконечном множестве $Y=\dbZ_p^*$, и $\Step(Y,R)$  --  $R$-модуль всех локально-постоянных функций на $Y$ (в данном случае все ступенчатые функции непрерывны!).

Напомним, что {\it распределение}  $\mu$ на $Y$
со значениями в  $\Mc$ это конечно-аддитивная функция на открытых подмножествах  $U\subset Y$:
$$
\mu\colon\,
\left\{
\begin{matrix}
\hbox{открытые подмножества }\cr
U\subset Y\cr
\end{matrix}
\right\}\longrightarrow\Mc.
$$
Другими словами, $\mu$ -- это гомоморфизм $R$--модулей
$$
\mu: \Step(Y,R)\to \Mc
$$
Напомним, что {\it мерой} на $Y$ со значениями в $\Mc$ называется {\it непрерывный}
 гомоморфизм $R$-модулей
$$ 
\mu : \Cr(Y,R) \longrightarrow \Mc.  
$$

Ограничение  $\mu$ на $R$--подмодуль  $\Step(Y,R) \subset \Cr(Y,R)$
определяет распределение, обозначаемое той же буквой $\mu$, причём  мера $\mu$ 
однозначно определена по соответствующему распределению, поскольку  $R$--подмодуль $\Step(Y,R)$ 
``плотен'' в $\Cr(Y,R)$. 
Это утверждение выражает  общий  факт
о равномерной непрерывности непрерывной функции на компакте $Y$.


\begin{corollary}[Мазур]\label{cor2}
Существует единственная $\dbQ_p$=значная мера  $\mu^{(c)}$
на $\dbZ_p^*$ такая, что для всех  $k\ge 1$ имеем 
$\int\limits_{\dbZ_p^*}x^kd\mu^{(c)}=
\zeta_{(p)}^{(c)}(1-k)=(1-c^k)
(1-p^{k-1})
\zeta(1-k)$. Заметим, что $ \zeta(0)=-\frac 1 2,$ но $\zeta_{(p)}^{(c)}(0)=0$.
\end{corollary}

Действительно, если интегрировать  $h(x)$ по $\ZZ_p^{\times}$,
 то получается сравнение Куммера из теоремы \ref{THkor1.2}. 
С другой стороны, для определения меры, удовлетворяюшей 
 условиям следствия, определим интеграл $\int_{\ZZ_p^{\times}} \phi(x) \mu^{(c)}$
 для каждой непрерывной функции  
$\phi : \ZZ_p^{\times} \rightarrow \ZZ_p$. 
Для этого используется приближение непрерывной функции 
$\phi(x)$ многочленами (для них интеграл задан по определению),
затем остаётся перейти к пределу.
{\it Сравнения Куммера} из теоремы \ref{THkor1.2}.
показывают, что предел корректно определён, и даёт интеграл для любой непрерывной функции.

\subsection{$p$-адическая дзета-функция Куботы -- Леопольдта}

\subsubsection{Область определения $p$-адических дзета-функций}
Областью определения  комплексных дзета функций является группа
$$
\CC =\Hom_{\cont}(\R^{\times }, \CC^{\times }), \ \ s\mapsto (y \mapsto y^s).
$$  
По аналогии с классическим комплексным случаем областью определения 
 $p$-адических дзета-функций является $p$-адическая группа
  $$
   X_p = \Hom_{\cont}(\ZZ_p^{\times }, \CC_p^{\times }),
  $$ 
состоящая из всех непрерывных гомоморфизмов
проконечной группы $\ZZ_p^{\times}$ в мультипликативную группу поля Тэйта,
$\CC_p = \widehat{\overline{\QQ}}_p$ 
(пополнение алгебраического замыкания поля  $p$-адических чисел $\QQ_p)$. 
Мы будем рассматривать целые числа $k$ как гомоморфизмы $x_p^k : y \mapsto y^k$.

Конструкция Куботы и Леорольдта даёт существование 
$p$-адической аналитической функции  $\zeta_p : X_p \rightarrow \CC_p$  с единственным простым полюсом в точке $x = x_p^{-1}$, 
которая становится ограниченной аналитической функцией на  $X_p$
 после умножения на регуляризующий множитель
 $(c x(c) - 1)$, $(x \in  X_p, c \in \ZZ_p^{\times}$);
эта функция однозначно определена условием
\begin{eqnarray}\label{zp}
 \zeta_p(x_p^k) = (1 - p^k) \zeta(-k) ~~~~ (k \geq 1).  
\end{eqnarray}

Этот результат имеет очень естественную интерпретацию в рамках теории 
$p$-адического интегрирования
(в стиле результата Мазура, см. следствие \ref{cor2})

Замечательное свойство этой конструкции состоит в том, что она применима и для всех
характеров Дирихле  $\chi$ по модулю степени простого числа $p$.
Зафиксируем вло\-же\-ние 
\begin{eqnarray}\label{ip}
i_p : \overline{\QQ} \hookrightarrow \CC_p
\end{eqnarray}
и будем отождествлять поле $\overline{\QQ}$ с подполем в $\CC$ и в $\CC_p$.
Тогда характер Дирихле вида  
$$
\chi : (\ZZ/{p^N}\ZZ)^{\times} \rightarrow \overline{\QQ}^{\times}
$$
становится элементом подгруппы кручения 
$$
X_p^{\tors} \subset X_p= \Hom_{\cont}(\ZZ_p^{\times }, \CC_p^{\times })  
$$
 и равенство (\ref{zp}) остаётся в силе и для специальных значений
  $L(-k,\chi )$ соответствующих 
$L$-рядов Дирихле
$$
 L(s,\chi) = \sum_{n=1}^{\infty} \chi(n) n^{-s}
               = \prod_{{\displaystyle \ell}  \small\hbox{ простые }\atop \hbox{ числа }} 
               (1-\chi(\ell) \ell^{-s})^{-1}, 
$$
при этом мы имеем 
 \begin{eqnarray}\label{lp}\hspace{0.3cm}
 \zeta_p(\chi x_p^k) = i_p \left[ (1 - \chi(p) p^k) L(-k,\chi) \right]
                           \  (k \geq 1,~~ k \in \ZZ, ~~
                           \chi \in X_p^{\tors}).  \hskip0.4cm
\end{eqnarray}
\subsubsection{Неархимедово преобразование Меллина}

Пусть $\mu$ обозначает  ${\CC}_p$-значную меру  $\ZZ_p^{\times}$. Тогда 
{\it неархимедово преобразование Меллина} 
меры $\mu$ определяется равенством 
  \begin{equation}  \label{EQMeltrans}
    L_{\mu}(x) = \mu(x) = \int_{\ZZ_p^{\times}} x \dd \mu, ~~~~ (x \in X_p), 
  \end{equation}
и представляет некоторую ограниченную ${\CC}_p$-аналитическую функцию
  \begin{equation}
    L_{\mu} : X_p \longrightarrow {\CC}_p.
  \end{equation}
  Действительно, ограниченность функции $L_{\mu}$ очевидна поскольку все характеры  $x \in X_p$ при\-ни\-мают значения в $\oc_p$ и $\mu$ также ограничена. Аналитичность $L_{\mu}$ вырыжает общее 
  свойство интеграла
(\ref{EQMeltrans}), поскольку он аналитически зависит от параметра
$x \in X_p$. 
Можно доказать (теорема Ивасава), что ограниченные ${\CC}_p$-аналитические функции 
взаимно-однозначно соот\-вет\-ству\-ют  
${\CC}_p$-значным мерам $\mu$ на  $\ZZ_p^{\times}$  посредством неархимедова преобразования Меллина.

\subsubsection{Пример: $p$-адическая
дзета-функция Куботы -- Леопольдта}

Для меры Мазура из следствия \ref{cor2}
функция на $X_p$
  \begin{equation}
    \z_p(x) = \left(1 - c^{-1} x(c)^{-1} \right)^{-1} L_{\mu^{(c)}}(x)
           ~~~~ (x \in X_p)
  \end{equation}
однозначно определена и голоморфна
за исключением  единственного простого полюса в точке $x = x_p^{-1}$, 
и становится ограниченной аналитической функцией на  $X_p$
 после умножения на регуляризующий множитель
 $(c x(c) - 1)$, $(x \in  X_p, c \in \ZZ_p^{\times}$);
эта функция однозначно определена условием
(\ref{zp}).


\subsection*{Признательность автора}
Искренне благодарю Эрнеста Борисовича Винберга  за 
приглашение подготовить статью  для журнала ``Математическое Просвещение'' 2008, 
посвящённого $p$-адическим числам и их при\-ло\-же\-ни\-ям.

\bibliographystyle{plain}

\end{document}